\documentclass[12pt]{amsart}
\usepackage{amsmath, amssymb, euscript, hyperref}
\usepackage{graphicx}

\newtheorem{thm}{Theorem}[section]
\newtheorem{lem}[thm]{Lemma}
\newtheorem{prop}[thm]{Proposition}
\theoremstyle{definition}
\newtheorem{defn}[thm]{Definition}

\newtheorem{rem}[thm]{Remark}

\newcommand \Cal {\mathcal }

\begin{document}

\bibliographystyle{amsplain}

\address{Azer Akhmedov, Department of Mathematics,
North Dakota State University,
Fargo, ND, 58102, USA}
\email{azer.akhmedov@ndsu.edu}

\begin{center} {\bf {\Large Cayley graphs with an infinite Heesch number}} \end{center}

\bigskip

\begin{center} {\bf Azer Akhmedov} \end{center}

\medskip

{\small Abstract: We construct a 2-generated group $\Gamma $ such that its Cayley graph possesses finite connected subsets with arbitrarily big finite Heesch number.}   

\vspace{1cm}

\section{INTRODUCTION}

 A {\em Heesch number} of a polygon $P$ is the maximum number layers of polygons isometric to $P$ that can surround $P$ without overlapping. For example, a rectangle in the plane has Heesch number infinity since it actually tiles the entire plane. On the other hand, there exist polygons which form only a {\em partial tile} for the plane, hence they have a finite Heesch number. It is more interesting to find polygons with a finite Heesch number so we will drop the perfect tiles of the plane from the consideration. 
 
 \medskip 
 
 The term Heesch number is named after the geometer Heinrich Heesch who found an example of a polygon with Heesch number 1 (See \cite{H}, p.23). This polygon is described in Figure 1 \footnote{Figure 1 and Figure 2 in this paper have been borrowed from the web page http://math.uttyler.edu/cmann/math/heesch/heesch.htm}, and consists of a union of a square, an equilateral triangle and a right triangle with angels 30-60-90. It is already much harder to find polygons with a Heesch number 2 (a well known example of such a polygon is given in Figure 2).
 
 \medskip
 
 The first examples of polygons with Heesch numbers 2, 3, 4 and 5 have been first discovered by A.Fontaine \cite{F}, R.Ammann, F.W.Marshall and C.Mann \cite{M1} respectively (See \cite{M2} for the history of this problem). 
  
\begin{figure}
  \includegraphics[width = 3in, height = 2in]{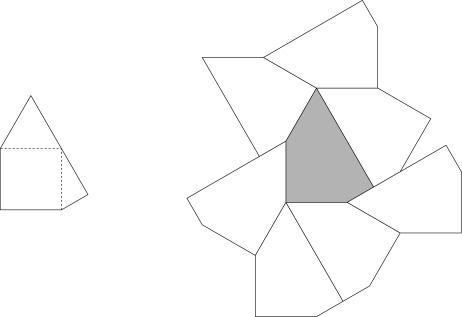}
\caption{A pentagon with Heesch number 1.}
\label{labelname}
\end{figure}

\medskip

 It is not known if there exists any polygon with a finite Heesch number $N > 5$. This is called {\em Heesch's Problem} and has close connections to a number of problems in combinatorial geometry such as {\em Domino Problem} and {\em Einstein Problem}. The latter asks if there exists a tile of the plane consisting of a single polygon $P$ such that any tiling of the plane by $P$ is aperiodic (``Ein  Stein" stands for ``one stone" in German; this word play is attributed to Ludwig Danzer). Such polygons have been constructed by G.Margulis and S.Mozes \cite {MM} in the hyperbolic plane, while in the Euclidean plane no such examples are known. (The famous Penrose tiling is known to be always aperiodic, but it uses two polygons, not one.) Incidentally, in the hyperbolic plane $\mathbb{H}^2$, the Heesch's Problem is solved completely; it is shown by A.S.Tarasaov \cite{T} that there exist polygons in $\mathbb{H}^2$ with an arbitrary Heesch number $N\geq 1$.    

\medskip

\begin{figure}
  \includegraphics[width = 5in, height = 3in]{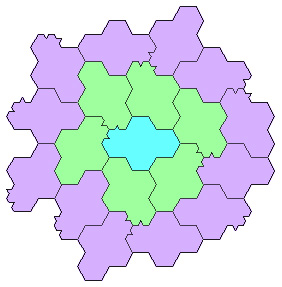}
\caption{C.Mann's example of a 20-gon with Heesch number 2.}
\label{labelname}
\end{figure}

\medskip

 In the Euclidean plane, one can also try to work with its lattice $\mathbb{Z}^2$ noticing that any finite connected set $K$ in the Cayley graph of $\mathbb{Z}^2$  with respect to the standard generating set $\{(\pm 1, 0), (0, \pm 1)\}$ gives rise to the polygon $$P(K) = \{(x,y)\in \mathbb{R}^2 \ | \ \displaystyle \min _{(u,v)\in K}\mathrm{max}\{|x-u|, |y-v|\} \leq \frac{1}{2}\}$$ which consists of the $\frac{1}{2}$-neighborhood of the discrete set $K$ in the $l^{\infty }$-metric. Then one can ask if there exist polygons of the form $P(K)$ with a big Heesch number. Indeed, the polygon with Heesch number 2 constructed in \cite{F} is of the form $P(K)$ but this is not the case for the known examples of polygons with Heesch numbers 3, 4 and 5.  Notice also that if $K\subset \mathbb{Z}^2$ is not connected then $P(K)$ is not a polygon any more, hence restricting to connected sets is natural.

 \medskip

 The notion of a partial tile can be defined in any group (where one moves the sets around by a left translation), and more generally in any graph (where one moves the sets around by an automorphism of the graph). It becomes interesting then if there exists a homogeneous graph (e.g. a vertex transitive graph), and more specifically, a Cayley graph which possesses connected partial tiles with arbitrarily big Heesch number. In the current paper we provide a positive answer to this question; we construct a group $\Gamma $ generated by a two element subset $S$ such that the Cayley graph of $\Gamma $ with respect to $S$ possesses partial tiles of arbitrarily big finite Heesch number (in the Cayley graph, the sets are moved around by left translations of the group). Notice that it is very easy to find disconnected sets with a big finite Heesch number, so without the connectedness condition the question is easy (and somewhat unnatural). 
 
 \medskip 
 
 To state our main result, we need to define the notions of tile, partial tile, and Heesch number in the setting of an arbitrary finitely generated group.
 
 \medskip  

   \begin{defn}(tiles) Let $\Gamma $ be a countable group, $F$ be a finite subset with cardinality at least 2. $F$ is called a tile if there exists $C\subset \Gamma $ such that $\Gamma = \displaystyle \mathop{\sqcup }_{g\in C}gF$. The partition $\displaystyle \mathop{\sqcup }_{g\in C}gF$ is called {\em a tiling}, and the set $C$ is called {\em the center set} of this tiling. We will always assume that $1\in C$. 
   \end{defn}
   
   \medskip
   
   We will also be interested in partial tiles of groups.
  
  \medskip
   
  \begin{defn}(partial tiles) Let $\Gamma $ be a countable group, $F$ be a finite subset with cardinality at least 2 such that $F$ does not tile $\Gamma $. For a subset $M$ of $\Gamma $, we say $F$ tiles $M$ if there exists $C\subset \Gamma $ such that $1\in C$ and $M \subseteq \displaystyle \mathop{\sqcup }_{g\in C}gF$. $F$ will be called {\em a partial tile} of $\Gamma $, and the partition $\displaystyle \mathop{\sqcup }_{g\in C}gF$ is called {\em a partial tiling}.
  \end{defn}
   
    \medskip
    
   Let now $\Gamma $ be a finitely generated group. We will fix a finite symmetric generating set $S$ of $\Gamma $, and study partial tiles in the Cayley graph of $\Gamma $ w.r.t. the left invariant Cayley metric given by $S$. For all $g\in \Gamma $, $|g|_{\Gamma }$ will denote the length of the element $g$ in the Cayley metric, and for all $x, y\in G$, $d_{\Gamma }(x,y)$ will denote the distance between $x$ and $y$, i.e. $d_{\Gamma }(x,y) = |x^{-1}y|$ (we will drop the index if it is clear from the context which group we are considering). For any $g\in \Gamma $, we will also write $B_g(r) = \{x\in \Gamma \ | \ d(g,x)\leq r\}$ for the {\em ball of radius $r$ around $g$}; for any two subsets $A, B$ of $\Gamma $ we will write $d(A,B) = \displaystyle \min _{x\in A, y\in B}d(x,y)$ for the {\em distance between the sets $A$ and $B$}; and for any finite subset $A\subset \Gamma $, we will also write $\partial A = \{x\in A \ | \ d(x,y) = 1 \ \mathrm{for \ some} \ y\in \Gamma \backslash A\}$ for the {\em boundary of $A$}. 
   
   \medskip 
  
   For a partial tile $K$, we say that the Heesch number of $K$ equals $N$, if one can tile $N$ layers around $K$ but not $N+1$ layers. To be precise, we need the following definitions.
   
   \medskip
   
   \begin{defn} (layers) Let $C\subseteq \Gamma $ such that $\pi = \displaystyle \mathop{\sqcup }_{g\in C}gK$ is a partial tiling. Let also $C_0\subseteq C$. We say $C_0K$ is the layer of level 0 of $\pi $ if $C_0 = \{1\}$. For a subset $C_1\subseteq C$, we say $C_1K$ is the layer of level 1 of $\pi $, if $C_1$ is a minimal subset of $C$ such that $C_1\cap C_0 = \emptyset $ and $\{x\in \Gamma   \ | \ d(x,K) = 1\}\subseteq C_1K$. (notice that if $C_1$ exists then it is unique).
   
   For any $n\geq 2$, inductively, we define the layer of level $n$ of $\pi $ as follows: if $C_0K, C_1K, \dots , C_{n-1}K$ are the layers of level $0, 1, \dots , n-1$ respectively, then we say $C_nK$ is the layer of level $n$ if $C_n$ is a minimal subset of $C$ such that $$C_n\cap \displaystyle \mathop{\sqcup }_{0\leq i\leq n-1}C_i = \emptyset \ \mathrm{and} \ \{x\in \Gamma   \ | \ d(x,\mathop{\sqcup }_{0\leq i\leq n-1}C_iK) = 1\}\subseteq C_nK$$
   \end{defn}
   
   \medskip
   
   Motivated by the above definitions, one naturally defines the notion of a Heesch number for partial tiles of the group $\Gamma $.
   
   \medskip
   
   \begin{defn} (Heesch number) Let $K$ be a finite subset of $\Gamma $. We write Heesch$(K) = N$ if $N$ is the maximal non-negative integer such that $\Gamma $ has a partial tiling by the left shifts of $K$ which has $N$ layers.
   \end{defn}     
   
    \medskip
    
  If a finite subset $K$ tiles $\Gamma $ then we will write Heesch$(K) = \infty $. 
  
  \medskip
  
  The main result of the paper is the following 
  
  \medskip
  
  \begin{thm} \label{thm:main}There exists a finitely generated group $\Gamma $ with a fixed finite generating set $S$ such that for any natural $N$, $\Gamma $ has a connected partial tile $K_N$ with a finite Heesch number bigger than $N$.   
   \end{thm}
    
    \medskip
  
  Let us emphasize again that without the condition ``connected" the result would be trivial as the group $\mathbb{Z}$ easily possesses disconnected partial tiles with arbitrarily big finite Heesch number. The statement of Theorem \ref{thm:main} also motivates the notion of a Heesch number for an arbitrary Cayley graph $\Cal{G} = \Cal{G}(\Gamma , S)$ of a group $\Gamma $ with respect to the symmetric generating set $S$, namely, the Heesch number of $\Cal{G}(\Gamma , S)$ is the maximal $N$ such that $\Cal{G}(\Gamma , S)$ possesses a finite connected partial tile $K_N$ with Heesch number $N$. Notice that this notion is quite sensitive to the choice of the generating set $S$. 
    
    \bigskip
    
    \section{Hyperbolic Limits} 
    
    We will be using the well known concept of hyperbolic limits. The reader may consult with \cite{GH} for basic notions of the theory of word hyperbolic groups but we will assume nothing other than the familiarity with the definition of a word hyperbolic group. Following the convention, we will say that a word hyperbolic group is {\em elementary} if it is virtually cyclic. Let us first recall a well known theorem due to Gromov and Delzant which motivates the notion of hyperbolic limit.  
    
    \medskip
    
    \begin{thm} \label{thm:gromov}(See \cite{G} and \cite{D}) Let $H$ be a non-elementary word hyperbolic group with a fixed finite generating set. Then for any non-torsion element $\gamma \in H$ and for any $R>0$ there exists a positive integer $N_0$ such that for all $N > N_0$ the quotient $H' = H/\langle \gamma ^N = 1 \rangle $ is non-elementary word hyperbolic, moreover the quotient map $\pi : H\rightarrow H'$ is injective on the ball $B_R(1)$ of radius $R$ around the identity element. [in other words, adding the relation $\gamma ^N = 1$ is injective on the ball of radius $R$ and the quotient remains non-elementary word hyperbolic]. 
    \end{thm}
    
    \medskip
    
    Let now $H$ be a non-elementary word hyperbolic group with a fixed finite symmetric generating set $S$. The group   $H_{\infty }$ is called a hyperbolic limit of $H$ if there exists a sequence $H_0, H_1, H_2, \ldots $ of non-elementary word hyperbolic groups such that 
    
    \medskip
    
    (i) $H_0 = H$;

  \medskip
  
    (ii) $H_{n+1}$ is a quotient of $H_n$ for all $n\in \mathbb{N}\cup \{0\}$;
  
  \medskip
    
    (iii) for all $n\in \mathbb{N}\cup \{0\}$, the quotient epimorphism $\pi _n : H_n\rightarrow H_{n+1}$ is injective on the ball of radius $n+1$ around identity element w.r.t. the generating set $S$ [more precisely, with respect to   the generating set $\pi _{n-1}\ldots \pi _1\pi _0(S)$, but by  abusing the notation, we will denote it by $S$].    
  
  \medskip
  
  Since the ball of radius $n$ remains injective (unchanged) by all the epimorphisms $\pi _i , i > n$, the union of these stable balls determines a group, denoted by $H_{\infty }$, and called {\em a hyperbolic limit} of $H$. If $g\in H_n$ belongs to the ball of radius $n$ around the identity element then, by abusing the notation, we will denote the image of $g$ in $H _{\infty }$ by $g$. So the image of the generating set $S$ of $H$ in $H_{\infty }$ will be denoted by $S$.
  
  \medskip
  
  A non-elementary word hyperbolic group may have many different hyperbolic limits, and groups which are very far from being word hyperbolic can be hyperbolic limits. Hyperbolic limits are very useful; for example, using Theorem \ref{thm:gromov}, one immediately obtains a finitely generated infinite torsion group as a hyperbolic limit of an arbitrary non-elementary word hyperbolic group. 
  
  \medskip
  
  In the proof of Theorem \ref{thm:main}, the group $\Gamma $ will be constructed as a hyperbolic limit of virtually free groups. Having a virtually free group at each step allows a great simplification in the argument but that also means we need to make extra efforts to keep the group virtually free at each step. Indeed, with a much more complicated argument, starting with an arbitrary non-elementary word hyperbolic group $H_0$, one can construct a group as a quotient of $H_0$ with connected partial tiles of arbitrarily big finite Heesch number.  
  
  \vspace{1cm}
  
   \section{Intermediate Results}
   
   \bigskip
   
   The following simple lemma will be extremely useful.   
   
   \begin{lem} \label{lem:quotient} Let $F$ be a finite subset of $\Gamma$, and $\pi :\Gamma \rightarrow \Gamma '$ be an epimorphism such that $\pi (F)$ is a tile of $\Gamma '$ and $|F| = |\pi (F)|$. Then $F$ is a tile of $\Gamma $.
   \end{lem}
   
   \medskip
   
   {\bf Proof.} Assume $N = ker(\pi )$, and $\Gamma ' = \displaystyle \mathop{\sqcup } _{g'\in C'\subseteq \Gamma '}g'\pi (F)$ is a partition of $\Gamma '$ into tiles. For every $g'\in \Gamma '$ we choose a representative $g\in \Gamma $ with $\pi (g) = g'$. Let $C$ be the set of all representatives. Then we have a partition $\Gamma = \displaystyle \mathop{\sqcup }_{g\in C, n\in N}gnF$. $\square $
   
   \bigskip
   
   In the proof of the main theorem, we will be considering virtually free groups. If $G$ is such a group with a fixed Cayley metric $|.|$, and $N$ is finite index free subgroup of rank $r\geq 2$ with a generating set $S$ of cardinality $r$, then, in general, it is possible that $|g| < |s|$ where $s\in S$ while $g\in N\backslash (\{1\}\cup S\cup S^{-1})$. However, we can avoid this situation by taking $N$ to be a very deep (still of finite index) subgroup in $G$. The following lemma is a simple consequence of hyperbolicity and residual finiteness of finitely generated virtually free groups.    
   
   \begin{lem}\label{lem:minimal} Let $G$ be a virtually non-abelian free group with a fixed finite generating set and the corresponding left-invariant Cayley metric $|.|$. Then there exists a finite index normal subgroup $N\unlhd G$ such that the following conditions hold:
   
   (c1) $N$ is a free group of rank $k\geq 2$ generated by a subset $\{g_1, \dots , g_k\}$, 
   
   (c2) $\min \{|g_1|, \dots , |g_k|\} = \mathrm{min}\{|x| \ : \ x\in N\backslash \{1\}\}$. $\square $ 
   \end{lem}
    
  \medskip
  
  {\bf Proof.} Since $G$ is virtually free, with respect to the metric $|.|$, it is $\delta $-hyperbolic for some $\delta > 0$. Let $G_0$ be a free subgroup of $G$ of finite index such that $G_0\cap B_{10\delta }(1) = \{1\}$. Let $r = \mathrm{rank}(G_0)$ and $S_0 = \{f_1, \dots , f_r\}$ be a generating set of $G_0$. We will denote the left invariant Cayley metric of $G_0$ with respect to $S_0$ by $|.|_0$.
  
  \medskip
  
  $G_0$ has a finite index subgroup $N$ such that $N$ is a normal subgroup of $G$. Let $k = \mathrm{rank}N$, and $$\mathcal{A} = \{(\gamma _1, \dots , \gamma _k) \ : \  \langle \gamma _1, \dots , \gamma _k \rangle = N, |\gamma _1|_0 \leq |\gamma _2|_0 \leq \dots \leq |\gamma _k|_0\}.$$ 
  
  Now we choose an arbitrary $(g_1, \dots , g_k)\in \mathcal{A}$ such that for any $(\gamma _1, \dots , \gamma _k)\in \mathcal{A}$ we have $|g_i|_0\leq |\gamma _i|_0, 1\leq i\leq k$. 
  
  \medskip
  
  Assuming the opposite, let $g\in N\backslash \{1\}$ such that $|g| = \min \{|x| \ : \ x\in N\backslash \{1\}\}$ and $g$ does not belong to any generating set of $N$ of cardinality $k$. We can write $g$ as a reduced word in the alphabet $\{g_1^{\pm 1}, \dots , g_k^{\pm 1}\}$. Let $p$ the minimal index such that $g_p$ or $g_p^{-1}$ occurs in $W$. Then, necessarily, $|g|_0 > |g_p|_0$. But this inequality contradicts $\delta $-hyperbolicity of $G$ with respect to the metric $|.|$. $\square $
  
  \medskip
  
  We now would like to state the central result of this section. 
  
  \medskip
  
  \begin{prop} \label{prop:bigger3} Let $G$ be a finitely generated virtually non-abelian free group, $N$ be a finite index normal subgroup of $G$ such that $N$ is a free group of rank $k\geq 4$ generated by elements $g_1, \dots , g_k$. Then there exists a generating $k$-tuple $(h_1, \dots , h_k)$ of $N$ such that for all positive real numbers $R>0$, there exists $n\geq 1$ such that the quotient $G /\langle h_i^n = 1, 1\leq i\leq k\rangle $ is a virtually non-abelian free group and the quotient epimorphism $G \to G /\langle h_i^n = 1, 1\leq i\leq k\rangle $ is injective on the ball of radius $R$ around the identity element.
  \end{prop}
  
  \medskip
  
  Before proving Proposition \ref{prop:bigger3}, we need to make a small digression into the outer automorphisms of free groups explaining also why do we need the condition $k\geq 4$. Let $G, N$ be as in Proposition \ref{prop:bigger3}, i.e. $G$ is a virtually non-abelian free group, and $N$ is a finite index free normal subgroup. If $g\in N$ then we have $\Cal{N} _N(g)\leq \Cal{N} _G(g)\leq N$ but the normal closure $\Cal{N} _G(g)$ of $g$ in $G$ can be much bigger than the normal closure $\Cal{N} _N(g)$ of $g$ in $N$. This is an undesirable situation for us; however, one can replace $(g_1, \dots , g_k)$ with another generating $k$-tuple $(h_1, \dots , h_k)$ of $N$, and let $g$ be quite special by taking $g = h_1^m$, for some (sufficiently big) $m\geq 1$. Then we can take the normal closure of a more symmetric set $S = \{h_1^m, \dots , h_k^m\}$, and it turns out that if $k\geq 4$ then the normal closure of $S$ in $N$ coincides with its normal closure in $G$. To see this let us first recall the following nice result of B. Zimmerman \cite{Z}.
  
  \medskip
  
  \begin{thm}\label{thm:zimmermann} Let $k\geq 2$. A finite subgroup $Out (\mathbb{F}_k)$ has a maximal order 12 for $k = 2$, and a maximal order $2^kk!$ for $k\geq 3$. Moreover, for $k\geq 4$, all finite subgroups are conjugate to a unique maximal subgroup $H_k$ of order $2^kk!$. $\square $
  \end{thm}      
  
  If $a_1, \dots , a_k$ are some generators of $\mathbb{F}_k$ then let $$L_k = \{\phi \in Aut \mathbb{F}_k \ | \ \forall i, \phi (a_i) \in \{a_1^{\pm 1}, \dots , a_k^{\pm 1}\}.$$ Notice that $L_k$ is a subgroup of $Aut \mathbb{F}_k$ of order exactly $2^kk!$. Moreover, the group $L_k$ induces a finite subgroup $\overline{L_k}$ of $Out (\mathbb{F}_k)$ of the same cardinality. Then by Theorem \ref{thm:zimmermann}, for $k\geq 4$, all finite subgroups of $Out (\mathbb{F}_k)$ are conjugate to a subgroup of $\overline{L_k}$ (in other words, in the statement of Theorem \ref{thm:zimmermann}, one can take $H_k$ to be $\overline{L_k}$). 
   
 \medskip
 
 Thus we obtain the following lemma.
 
 \begin{lem} \label{lem:luiza} $G$ be a finitely generated virtually non-abelian free group, $N$ be a finite index normal subgroup of $G$ such that $N$ is free of rank $k\geq 4$. Then there exists a generating $k$-tuple $(h_1, \dots , h_k)$ of $N$ such that for all  $m\geq 2$, we have $\Cal{N} _N(S) = \Cal{N} _G(S)$  where $S = \{h_1^m, \dots , h_k^m\}$.
 \end{lem}
 
 \medskip
 
 The Lemma \ref{lem:luiza} guarantees that $N/\Cal{N} _N(S)$ is a finite index subgroup of $G/\Cal{N} _G(S)$; it remains to recall a folklore result that $N/\Cal{N} _N(S)$ is always virtually free. For the sake of completeness we would like to formalize this claim in the following lemma. 
 
 \medskip
 
 \begin{lem}\label{lem:involution} Let $N\cong \mathbb{F}_k$ be a free group of rank $k\geq 3$ generated by elements $a_1, \dots , a_k$. Then for all $m\geq 2$, the quotient group $N_m = N/\langle a_1^m = \dots = a_k^m = 1\rangle $ is virtually non-abelian free.
 \end{lem}   
 
 \medskip
 
 For the proof, it suffices to recall a well known fact that the free product of two finitely generated virtually free groups is still virtually free (in particular, the free product of finite groups is virtually free), moreover, the free product of two non-trivial finite groups is virtually non-abelian free provided at least one of these finite groups have order at least 3.  
 
 \medskip   
 
 Now, Proposition \ref{prop:bigger3} follows immediately from Theorem \ref{thm:gromov}, Lemma \ref{lem:luiza} and Lemma \ref{lem:involution}, by taking $m$ to be a sufficiently big even number.
 
 \medskip
 
  We close this section with the following very useful result.
  
  \begin{lem}\label{lem:connected} Let $G$ be a group with a fixed finite symmetric generating set $S$ of cardinality at least four such that $S$ has no relation of length less than four, and for some $s\geq 1$, and for all distinct $x, y, z\in S$ there exists a path in $G\backslash \{1\}$ of length at most $s$ from $x$ to $y$ not containing $z$. Let also  $N$ be a finite index normal  subgroup of $G$ such that $N\cap B_{10s}(1) = \{1\}$, and $g\in N$ such that $|g| = \mathrm{min}\{|x| \ | \ x\in N\backslash \{1\}\}$. Then there exists a connected subset $A\subset G$ such that the following conditions hold:
  
  (a1) $|A| = |G/N|$, 
  
  (a2) $1\in A$, 
   
  (a3) $x^{-1}y\notin N$ for all distinct $x, y\in A$,
  
  (a4) $d(g,A) = 1$,
   
  (a5) $\pi = \displaystyle \mathop{\sqcup }_{h\in N}hA$ is the only left tiling of $G$ by $A$.
    
  \end{lem} 
  
  {\bf Proof.} Let $\epsilon : G\rightarrow G/N$ be the quotient map, and $\mathcal{G}, \mathcal{G}_1$ be the Cayley graphs of the groups $G, G/N$ with respect to the generating sets $S, \epsilon (S)$ respectively. Let also $r = (x_0 = 1, x_1, \dots , x_n, x_{n+1} = g)$ be a path in $\mathcal{G}$ connecting $1$ to $g$ such that $|x_i| = i, 0\leq i\leq n+1$. 
  
  \medskip
  
  We will first demonstrate the construction of the set $A$ satisfying conditions (a1)-(a4). In the Cayley graph $\mathcal{G}_1$, we consider the path $r_1 = (1, \epsilon (x_1), \dots , \epsilon (x_n))$, and build the subsets $B_1, B_2, \dots , B_{|G/N|-n}$ in $G/N$ inductively as follows.
  
  \medskip
  
  We let $B_1 = \{1, \epsilon (x_1), \dots , \epsilon (x_n)\}$, and if the sets $B_1, \dots , B_k$ are already defined for some $k < |G/N| - n$, then we let $B_{k+1} = B_k\sqcup \{z_k\}$ where $|z_k^{-1}z| = 1$ for some $z\in B_k$ (i.e. the distance from $z_k$ to $B_k$ equals 1, in the Cayley graph $\mathcal{G}_1$). Then $B_{|G/N|-n} = G/N$, and we start defining $A_1, \dots , A_{|G/N|-n}$ inductively as follows: we let $A_1 = \{1, x_1, \dots , x_n\}$, and if $A_1, \dots , A_k$ are defined for some $k < |G/N| - n$, then we let $A_{k+1}$ to be any connected set in $\mathcal{G}$ such that $A_{k+1} = A_k\cup \{y_k\}$ where $y_k\in \epsilon ^{-1}(z_k)$. 
  
  \medskip
  
  Then the set $A_{|G/N|-n}$ satisfies conditions (a1)-(a4) of the lemma. To make it satisfy the condition (a5) as well we need to modify our strategy little bit. 
  
  \medskip
  
  Let $$a = x_n^{-1}x_{n+1} = x_n^{-1}g, b = x_0^{-1}x_1 = x_1, U_1 = \{x\in G \ | \  d(x,x_n) = 1, x\neq g\}.$$ Let also $U_2\subset B_2(x_n)\backslash B_1(x_n)$ such that for all $x\in U_1$, the set $$(B_1(x)\backslash \{x\})\cap U_2$$ consists of a single element $y(x)$ where $x^{-1}y(x)\notin \{a, a^{-1}\}$ (See Fig. 3). Finally, we let $U = U_1\sqcup U_2$ and $V_0 = \{1, x_1, \dots , x_{n-1}\}\sqcup U = A_1\sqcup U$. Notice that by minimality assumption on $|g|$, we have $b\neq a^{-1}$, $|V_0| < |G/N|$ and the map $\epsilon $ is still injective on the set $V_0$; so the set $V_0$ already satisfies conditions (a2),(a3) and (a4). Our goal is to extend $V_0$ such that it also satisfies (a1) and (a5). However, notice that $V_0$ is not connected so first we would like to make it connected. 
  
  \medskip
  
  \begin{figure}
  \includegraphics[width = 3in, height = 2in]{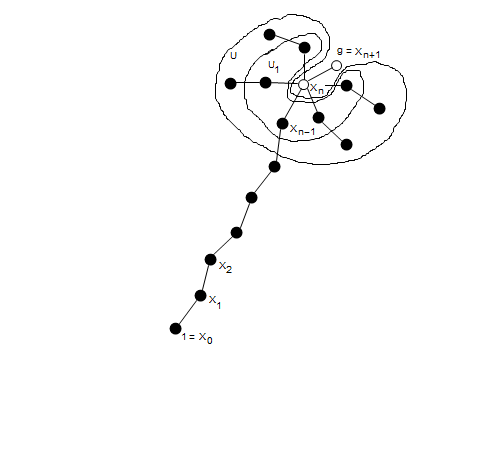}
\caption{The sets $V_0$ and $U$; the elements of $V_0$ are represented with black dots}
\label{labelname}
\end{figure}
  
   For this purpose, let $S\backslash \{a\} = \{g_1, \dots , g_{|S|-1}\}$, and for all $1\leq j\leq |S|-1$, let $R_j$ be a shortest path connecting $x_ng_j$ to $x_{n-1}$ avoiding $\{x_n,g\}$. Then we let $V = V_0\cup \displaystyle \mathop{\cup }_{1\leq j\leq |S|-1}R_j$, and observe that $V$ still satisfies conditions (a2), (a3) and (a4)\footnote{for the condition (a3), it sufficies to notice that $N\cap B_{10s}(1) = \{1\}$, and recall the condition on the generating set $S$: $S$ has cardinality at least four;  there is no relation of length less than four among elements of $S$; and for all distinct $x, y, z\in S$ there exists a path in $G\backslash \{1\}$ of length at most $s$ from $x$ to $y$ not containing $z$.}, moreover, $V$ is connected and $V\cap \{x_n,g\} = \emptyset $.  
    
    \medskip

  We let $A_1' = V, B_1' = \epsilon (V)$. Now for all $1\leq k < |G/N|-|V|$ suppose the connected subsets $A_1', \dots , A_k'$ of $G\backslash \{x_n, g\}$ are already defined such that $\{x\in \partial A_j' \ | \ xa\in A_j'\}\backslash \{a^{-1}\} = \emptyset $ for all $1\leq j\leq k$. Then we let $B_k' = \epsilon (A_k), 1\leq k < |G/N|-|V|$, and build a connected subset $A_{k+1}'\supset A_k'$ such that the conditions 
  
  (i) $1\leq |A_{k+1}'\backslash A_{k}'|\leq 2$,  
  
  (ii) $\{x\in \partial A_{k+1}' \ | \ xa\in A_{k+1}'\}\backslash \{a^{-1}\} = \emptyset $, 
   
  (iii) $A_{k+1}'\cap \{x_n, g\} = \emptyset $,
   
   hold.
   
\medskip
   
   For this purpose, let $$D = \{x\in G \ | \ d(x, A_k') = 1\}, D' = \{x\in G/N \ | \ d(x, B_k') = 1\}$$ and for all $z\in D'$, define $C(z) = \{x\in B_k' \ | \ d(x,z) = 1\}$. If there exist $z\in D', y\in D\backslash \{x_{n+1}, g\}, u\in C(z)$ such that $\epsilon (y) = z$ and $z\neq u\epsilon (a^{-1})$, then we define $B_{k+1}' = B_k\cup \{z\}$ and let $A_{k+1}' = A_k'\sqcup \{y\}$. 
  
 \medskip 
  
  But if such $z, y$ and $u$ do not exist, then, necessarily, there exist $z_1, z_2\in (G/N)\backslash B_k', y_1, y_2\in G\backslash (A_k'\sqcup \{x_n, g\})$ such that $\epsilon (y_i) = z_i, 1\leq i\leq 2, z_1\in D', d(z_2, z_1) = 1$ and $z_2^{-1}z_1\neq \epsilon (a^{-1})$. Then we let $B_{k+1}' = B_k'\sqcup \{z_1, z_2\}$, and define $A_{k+1}' = A_{k}'\sqcup \{y_1, y_2\}$. 
  
  \medskip
  
  Finally, let $m$ be such that $B_m' = G/N$. Then the set $A = A_m'$ satisfies conditions (a1)-(a5). $\square $
  
  \bigskip
 
 \begin{rem} Let us emphasize that because of a particular shape of the set $V_0$, if it tiles the group $G$ then we have a forced unique partial tiling in a large part of the group. This is because the set $V_0$ consists of a ``head" $U$ and a ``tail" $V_0\backslash U$; we have a small ``hole" $\{x_n, g\}$ at the head so that in any tiling of $G$ by $V_0$ the tail of a shift of $V_0$ must enter into this hole.  We use this observation as a key tool in the proof of Lemma \ref{lem:connected}. For the proof of the lemma, we need to extend the set $V_0$ (the extension is needed to satisfy the condition (a1)) by preserving this distinctive property of it which forces the uniqueness of the tiling. The existence of the tiling follows simply from the conditions (a1) and (a3).
 \end{rem}

 \section{Proof of the Main Theorem}
 
 Let $\mathbb{F}_2$ be generated by the set $\{a,b\}$, and $H$ be a virtually non-abelian free quotient of $\mathbb{F}_2$ such that there is no relation of length less than four among $a$ and $b$, moreover, for any distinct $x, y, z\in \{a,a^{-1},b,b^{-1}\}$ there exists a path in $H\backslash \{1\}$ connecting $x$ to $y$ and not containing $z$. Let also $s$ be the maximal length of all such paths. 
 
 \medskip
 
 For the proof, we will construct the group $\Gamma $ as a hyperbolic limit of the group $H$ with a fixed generating set $S = \{a,a^{-1},b,b^{-1}\}$. We will build the hyperbolic limit of $H = H_0$ inductively as follows. Suppose the groups $H_0, \dots , H_n$ have been constructed such that the following conditions hold:
 
 (i) $H_{i}$ is a quotient of $H_{i-1}$, for all $1\leq i\leq n$;
  
  \medskip
  
  (ii) $H_i$ is virtually non-abelian free, for all $0\leq i\leq n$;
  
  \medskip
    
  (iii) $H_i$ possesses a partial tile $K_i$ of finite Heesch number at least $i+1$, for all $0\leq i\leq n$; 
  
  \medskip
  
  (iv) For all $0\leq i\leq n-1$, the quotient epimorphism $\pi _i : H_{i}\rightarrow H_{i+1}$ is injective on the ball of radius $(i+1)(r(i)+1)$ around the identity element w.r.t. the generating set $S$ where $r(i) = \mathrm{max}\{|g|_j \ | \ 1\leq j\leq i, g\in K_j\}, \forall i\geq 1$ and $|.|_i$ denotes the left-invariant Cayley metric in $H_{i}$ with respect to the generating set $S$ (we let $r(0) = 10s$). 
  
  \medskip
 
 Now, by Lemma \ref{lem:minimal} there exists a normal subgroup $N\unlhd H_n$ of a finite index such that $N$ is a free group of rank $k\geq 4$ generated by elements $a_1, \dots , a_k$ and for all $g\in N$, if $|g|_n = \mathrm{min}\{|x|_n \ | \ x\in N\backslash \{1\}\}$ then $|g|_n = \min \{|a_1|_n, \dots , |a_k|_n\}$ (i.e. conditions (c1) and (c2) of Lemma \ref{lem:minimal} hold). Without loss of generality we may assume that $$|a_1|_n = \mathrm{min}\{|x|_n \ | \ x\in N\backslash \{1\}\}.$$ Then, by Lemma \ref{lem:connected}, there exists a connected set $A\subset H_n$ satisfying conditions (a1)-(a5). (Then, in particular, $d(a_1, u) = 1$ for some $u\in A$).  
 
 \medskip
 
 Let also $R_1 = \mathrm{max}\{|g|_n \ | \ g\in K_i, 1\leq i\leq n\}$ (we let $R_1 = 0$ if $n=0$) and $R_2 = \mathrm{max}\{|a_1^j|_n \ | \ 1\leq j\leq 10(n+1)\}$.
 
 \medskip
 
 By Proposition \ref{prop:bigger3} there exists a generating set $\{h_1,\dots , h_k\}$ of $N$, and an odd number $p_n\geq 1$ such that the quotient $$H_{n+1}:= H_n /\langle h_i^{p_n} = 1, 1\leq i\leq k\rangle $$ is a virtually non-abelian free group, and the quotient epimorphism $$H_n \to H_n /\langle h_i^{p_n} = 1, 1\leq i\leq k\rangle $$ is injective on the ball of radius $R := \mathrm{max}\{R_1^n, R_2\}$ around the identity element. Then notice that, by Lemma \ref{lem:quotient}, all the partial tiles $K_1, \dots , K_n$ inject into $H_{n+1}$, moreover, the images of $K_i$ have a finite Heesch number at least $i+1$ in $H_{n+1}$. Then we take $K_{n+1} = A\cup a_1A$. Notice that $K_{n+1}$ has a finite Heesch number at least $n+1$. Thus we can continue the process.
 
 \medskip
 
 It remains to notice that the Heesch number of $K_{n+1}$ in $H_{n+1}$ is at most $\frac{p_n}{2}$ thus it is finite. Then the hyperbolic limit group $H_{\infty }$ will be a group with desired properties, i.e. the connected sets $K_1, K_2, \dots $ will be all partial tiles with finite Heesch numbers $h_1, h_2, \dots $ such that $h_n\geq n$ for all $n\geq 1$. $\square $ 
  
 \bigskip
 
 {\em Acknowledgment:} I am thankful to I.Pak for encouraging me to write this paper, and to C.Mann for the figures in his website (Figure 1, 2) that I have borrowed.


\begin{thebibliography}{99}


\bibitem{D} Delzant, T. Sous-groupes distingues et quotients des groupes hyperboliques.
{\em Duke Math. J.}, {\bf 83}, (1996), no.3., 661-682.

\bibitem{F} Fontaine, A. \ An infinite number of plane figures with heesch number two, J. Combin. Theory Ser. A 57
(1991) 151–156.
 
\bibitem{G} Gromov, M. \ Essays in Group Theory, Math. Sci. Research Institue Publications,
8. edited by S.M.Gersten.

\bibitem{GH} Ghys, E. and de la Harpe, P. \ Sur les groupes hyperboliques d'apr\'es Mikhael Gromov. \ Birkh\"auser, 1990.

\bibitem{H} Heesch, H. \ Regulares Parkettierungsproblem, Westdeutscher Verlag, Cologne, 1968.

\bibitem{M1} Mann, C. \ Heesch's Problem and Other Tiling Problems, Ph.D. dissertation, University of Arkansas at
Fayetteville, 2001.

\bibitem{M2} Mann, C. \ Heesch's Tiling Problem \ {\em American Mathematical Monthly}, {\bf 111} (6), (2004) 509-517, 

\bibitem{MM} Margulis, G and Mozes, S. \ Aperiodic tilings of the hyperbolic plane by
convex polygons. {\em Israel J. of Mathematics} {\bf 107}, (1998), 319-325.

\bibitem{T} Tarasov, A. S. \ On the Heesch Number for the Hyperbolic Plane. \ {\em Mathematical Notes,} {\bf 88}, no.1, (2010) 97-104.

\bibitem{Z} Zimmermann, B. \ Finite groups of outer automorphisms of free groups. \ {\em Glasgow Math. J.} {\bf 38}, (1996) 275-282.
 

\end{thebibliography}
 \end{document}